\documentclass{amsart}

\usepackage{amsmath}
\usepackage{amsfonts}
\usepackage{amssymb}
\usepackage{amsthm}
\usepackage{hyperref}

\def\cal{\mathcal}

\newtheorem{thm}{Theorem}
\newtheorem{cor}{Corollary}
\newtheorem{lem}{Lemma}

\begin{document}

\title{Complexity results for CR mappings between spheres}

\author{John P.\@ D'Angelo}

\address{Dept.\@ of Mathematics, University of Illinois,
1409 W. Green St., Urbana, IL 61801, USA}
\email{jpda@math.uiuc.edu}

\author{Ji\v{r}\'i Lebl}

\address{Dept.\@ of Mathematics, University of Illinois,
1409 W. Green St., Urbana, IL 61801, USA}
\email{jlebl@math.uiuc.edu}

\maketitle

\begin{abstract}
Using elementary number theory, we
prove several results about the complexity of CR mappings between spheres.
It is known that CR mappings between spheres,
invariant under finite groups, lead to sharp bounds for 
degree estimates on real polynomials constant on a hyperplane.
We show here that there are infinitely many degrees for which the uniqueness of sharp examples fails.
The proof uses a Pell equation. We then sharpen our results and obtain various congruences guaranteeing nonuniqueness.
We also show that a gap phenomenon for proper mappings between balls does not occur beyond a certain target dimension.
This proof uses the solution of the postage stamp problem.

\medskip


\medskip

\end{abstract}

\section{Introduction}

There has been some recent interest in the complexity theory of CR mappings;
this paper deals with two specific questions in that area. In the first
we connect a uniqueness question about CR mappings between spheres to some 
combinatorial number theory. In the second we show that the {\it gap
phenomena} \cite{HJX2}
for proper mappings between balls applies only in the low codimension case, by establishing
that such gaps do not occur beyond a certain target dimension. Both results rely on 
elementary number theory; the first result uses the Pell equation and the second result uses 
Sylvester's formula for the solution of the postage stamp problem.

These results about CR mappings have interpretations 
in real algebraic geometry, and we continue the introduction
by discussing this connection. As usual let ${\bf R}[x]$ denote the polynomial ring 
in $n$ real variables, and let
${\bf R}_+[x]$ denote the cone of polynomials whose coefficients are non-negative.
We assume that $n\ge 2$ in order that the problem be interesting. Because
of the connection with CR mappings to spheres, elaborated in Section 2, we
are interested in the collection ${\cal H}(n)$ of elements of ${\bf R}_+[x]$
which take the constant value $1$ on the hyperplane defined by $\sum_{j=1}^n x_j =1$.

The following sharp result was proved in \cite{DKR}.
Suppose $p \in {\cal H}(2)$. Let $N$ denote the number
of distinct monomials in $p$ and let $d$ denote the degree of $p$.
Then $d \le 2N-3$ and for each odd $d=2r+1$ there is a polynomial of degree
$d$ with $2N-3=2(r+2)-3= 2r+1$.  One family of polynomials exhibiting this sharp bound 
has many interesting properties; we discuss some of these properties in Section 2.
They arose from considering CR mappings invariant under certain finite unitary groups.
See also \cite{D}, \cite{D3}, \cite{D4}, \cite{F2}, and \cite{LWW}.

We naturally ask whether the {\it sharp} polynomials are unique. The answer is no in general.
In Sections 3 and 4 we prove that non-uniqueness is a generic phenomenon.
In particular there are infinitely many odd $d$ for which
other inequivalent polynomials also exhibit the sharp bound. For $d=1$ and $d=3$
there is only one {\it sharp} polynomial; for $d=5$ there are two, 
but they are equivalent in a sense to be described. Theorem 2 tells us that for
each $d$ at least $7$ and congruent to $3$ mod $4$, 
there are inequivalent sharp polynomials. We also obtain in Corollary 2
a similar result when $d$ is congruent to $1$ mod $6$.
We say for simplicity that {\it uniqueness fails} in these cases.

Theorem 1 produces a general method for finding infinite sets of odd numbers for which 
uniqueness fails. The proof relies on the
integer solutions $(a,b)$ to the Pell equation $a^2 - 12 b^2 = 1$.
In particular we find infinitely many odd $d$ congruent to $1$ mod $4$
for which uniqueness fails. Theorem 2 gives precise formulas for additional
inequivalent sharp polynomials, but it requires 
considerable work to verify that their coefficients are positive. 
Corollary 2 gathers all these results.

The case of even degree is almost trivial; uniqueness
fails for every even degree at least $2$. See Section 4.
On the other hand we establish there
a complexity result in the case of even degree; the number of inequivalent
sharp polynomials (in two dimensions)
for a fixed even degree $2k$ tends to infinity as $k$ does.

In Theorem 3 from Section 5 we determine an integer $T(n)$ with the following property.
For each $N$ with $N\ge T(n)$ we can find a proper polynomial mapping from the unit ball
${\bf B}_n$ to ${\bf B}_N$ that cannot be mapped into a lower dimensional ball. This result shows that
the gap phenomenon goes away once the target dimension is sufficiently large. The proof is an elementary
construction using Sylvester's solution of the postage stamp problem.
The gap phenomenon does exist in low codimension. See \cite{Fa1} and \cite{HJX2} and Section 5.

It is natural to wonder, given the prominent role of the group invariant examples, whether any of our results 
are connected with the topology of Lens spaces. We do not consider this aspect of the problem in this paper.

\section{CR Mappings Between Spheres}

Let $S^{2n-1}$ denote the unit sphere in complex Euclidean space ${\bf C}^n$.
We will assume that $n\ge2$. A CR mapping between CR manifolds
is a continuously differentiable mapping that satisfies 
the tangential Cauchy-Riemann equations. See \cite{BER,D} for example.
For our purposes we may assume that a smooth CR mapping between spheres is the
restriction of a rational mapping (defined and holomorphic near the sphere).
We can do so because a non-constant CR mapping between spheres can be extended 
to be a proper holomorphic mapping between balls, and  a theorem of Forstneric
\cite{F1} 
guarantees that a proper holomorphic mapping between balls with sufficiently 
many derivatives at the boundary sphere is in fact a rational mapping.

For dimensions $n$ and $N$ with $2 \le n$
we denote by $P(n,N)$ the collection of rational mappings $f$ such that
$f: S^{2n-1} \to S^{2N-1}$. There is a natural equivalence relation in this
setting. Elements $f$ and $g$ in $P(n,N)$ are {\it spherically equivalent}
if there are automorphisms $\phi$ and $\psi$ of the unit balls in the domain and
target such that
$$ g = \psi \circ f \circ \phi. $$
When $f$ and $g$ are polynomial mappings in $P(n,N)$ and each preserves the origin,
then they are spherically equivalent only when they are unitarily equivalent.

For a fixed $n$, the set $P(n,N)$ grows increasingly
complicated as $N$ tends to infinity. The following
result \cite{D2} illustrates this complexity phenomenon. Given a polynomial $q:{\bf C}^n
\to {\bf C}$ such that $q\ne 0$ on $S^{2n-1}$,
there is an integer $N$ and a polynomial mapping 
$p:{\bf C}^n \to {\bf C}^N$ such that ${p \over q}$ is reduced to lowest terms and
such that ${p \over q}(S^{2n-1}) \subset S^{2N-1}$. Thus
${p \over q} \in P(n,N)$. Hence there is no restriction
on the possible denominators $q$ that can occur for elements of $P(n,N)$ if we choose
$N$ large enough. The proof of this result involves
ideas related to a complex variables analogue of Hilbert's seventeenth problem.
We also note the following easier fact. Fix $n$; given an integer $d$,
there exists an integer $N_0=N_0(n,d)$ such that, for $N \ge N_0$, $P(n,N)$ 
contains polynomials of degree $d$. The basic questions concern estimates on $d$ for elements of $P(n,N)$.

The following results are known; see \cite{DLP} for additional discussion and references.
Recall that we are assuming that $n\ge 2$. 

1) For $N < n$, all elements of $P(n,N)$ are constant and hence of degree zero.

2) \cite{P}. For $N=n$, all elements of $P(n,N)$ are constant or are linear fractional
transformations. Hence the degree is at most one.

3) \cite{Fa1}. For $n\le N \le 2n-2$, the degree of each element in $P(n,N)$ is at most one.

4) \cite{HJ}, \cite{HJX2}. For $3 \le n \le N \le 2n-1$, or for $4 \le n\le N \le 3n-4$, the degree of each element 
in $P(n,N)$ is at most two. 

5) \cite{Fa2}. For $N \le 3$, the degree of each element in $P(2,N)$ is at most three; 
there are four spherically inequivalent non-constant examples.

6) \cite{M}. The degree of each element in $P(2,N)$ is at most ${N(N-1) \over 2}$.

7) \cite{DKR}, \cite{D1}. For $n=2$ and $N$ at least three,
there is an element of $P(2,N)$ of degree $2N-3$. For $n\ge 3$
there is an element of $P(n,N)$ of degree $d$ with $d(n-1)= N-1$.  

8) \cite{DKR}. If $n=2$ and we restrict to monomial mappings in 7), 
then $2N-3$ is sharp.  Thus if $m \in P(2,N)$ is a 
monomial mapping, then its degree is at most $2N-3$.

9) \cite{DLP}. For $n\ge 2$ and monomial mappings in $P(n,N)$, (1) holds for the degree $d$:

$$ d \le {2n(2N-3) \over 3n^2-3n-2} \le {4 \over 3} {2N-3 \over 2n-3}. 
\eqno (1) $$

10) \cite{DLP}. For $n$ sufficiently large compared with $d$ and all monomials in $P(n,N)$,
the sharp inequality 2) holds for the degree $d$: 
$$ d \le {N-1 \over n-1}. \eqno (2) $$

We next recall (\cite{D}, \cite{DKR}, \cite{DLP}) how the restriction to monomial mappings 
allows us to express these questions
in terms of real polynomials with nonnegative coefficients.
Assume that $f:{\bf C}^n \to {\bf C}^N$ is a monomial mapping; thus
each component function of $f$ is a monomial, say $c_\alpha z^\alpha$ in
multi-index notation. If also
$f(S^{2n-1}) \subset S^{2N-1}$, then 
$$ \sum_\alpha |c_\alpha|^2 |z^\alpha|^2 = 1 $$ 
whenever $\sum_{j=1}^n |z_j|^2 =1$. Replacing $|z_j|^2$ by the real variable $x_j$
we obtain the equation
$$ \sum_\alpha |c_\alpha|^2 x^\alpha = 1 \eqno (3) $$ 
on the hyperplane $H$ given by $\sum_{j=1}^n x_j =1$.

Define $m$ by $m(x) = \sum_\alpha |c_\alpha|^2 x^\alpha$. Then $m$ is
a real polynomial with nonnegative coefficients and $m=1$ on the hyperplane
$H$; thus $m \in {\cal H}(n)$. Let ${\cal H}(n,d)$ denote the polynomials of degree $d$ in ${\cal H}(n)$.

Thus a monomial mapping in $P(n,N)$ of degree $d$ gives rise to 
an element of ${\cal H}(n,d)$ with $N$ distinct monomials.
Conversely, given $p \in {\cal H}(n,d)$ with $N$ distinct monomials,
we can reverse the procedure to find a monomial mapping of degree $d$ in $P(n,N)$.
The results 8) through 10) above therefore provide lower bounds
for the number of distinct monomials $N(p)$ in terms of the dimension
and degree of $p \in {\cal H}(n,d)$:

8') For $p \in {\cal H}(2,d)$, the sharp estimate $d \le 2N-3$ holds.

9') For all $n\ge 2$ and $p \in {\cal H}(n,d)$, estimate (1) holds.

10') For $n$ sufficiently large compared with $d$ and
$p \in {\cal H}(n,d)$, estimate (2) holds.

A family of polynomials exhibiting the sharp bound in 8') appear
in \cite{D}, \cite{D1}, \cite{DKR}, and \cite{DLP}. This particular family of polynomials has connections with
many branches of mathematics. For example, 
the polynomials are invariant under certain representations of 
finite cyclic groups (\cite{D1}, \cite{D3}), 
they have integer coefficients which provide a primality test \cite{D4}, they are closely
related to Chebychev polynomials, they arise in
problems such as de-nesting radicals in computer science, and so on. 

Given all these properties, it is not so surprising
that these polynomials provide the sharp bound in 8'). It is more
surprising that, for most odd $d$, there exist additional examples
exhibiting the sharp bound. Theorems 1 and 2 show how to construct
such examples of non-uniqueness. 

We now state precisely what we mean by {\it uniqueness} for degree $d$ and {\it equivalence} in ${\cal H}(2,d)$.
We say that $f$ and $g$ are equivalent in ${\cal H}(2,d)$ if either $f=g$ or $f(x,y)=g(y,x)$.
We say that {\it uniqueness holds} for a degree $d$ if the following is true.
Suppose $f,g \in {\cal H}(2)$ have the same degree and the same
number of terms. Furthermore, suppose that there is no $h \in {\cal H}(2)$ of the same degree
with fewer terms. Then $f$ and $g$ are equivalent in ${\cal H}(2,d)$.
Thus, uniqueness fails if we can find at least three different {\it sharp} polynomials in ${\cal H}(2)$.

\section{The Pell  Equation}

The reader should consult \cite{Len} for additional references and considerable discussion.
We discuss only a few issues considering the following Diophantine equation.
Let $ \lambda$ be a positive integer, assumed not to be a square. We
seek a pair $(d,k)$ of positive integers satisfying:
$$ d^2 = \lambda k^2 + 1 \eqno (4) $$

Given positive integers $d_1,k_1$ and a nonzero positive number $\lambda$ we write 
$r = d_1 + \sqrt{\lambda} k_1$ and ${\overline r}$ for the conjugate expression
${\overline r} = d_1 - \sqrt{\lambda} k_1$.
Then $d_1 = {r+{\overline r} \over 2}$ and 
$k_1 = { r - {\overline r} \over2\sqrt{\lambda}}$.

We define sequences of integers $d_m$ and $k_m$ by writing 
$$ (d_1 + \sqrt{\lambda} k_1)^m = d_m + \sqrt{\lambda} k_m. \eqno (5) $$
By standard methods for solving recurrences we obtain formulas for these integers:
$$ d_m = {r^m + {\overline r}^m \over 2} $$
$$ k_m = {r^m - {\overline r}^m \over 2 \sqrt{\lambda}}. \eqno (6) $$

We easily verify that if $(d_1,k_1)$ satisfies (4), then so does $(d_m,k_m)$.
In particular if we can find one solution we can find infinitely many.
We will need the following specific instance. For $\lambda =12$, note that $(7,2)$ is the solution
of (4) with $d_1$ and $k_1$ minimal. 

\begin{lem}
Put 
$$ d_m = {(7 + \sqrt{48})^m + (7 - \sqrt{48})^m \over 2}. \eqno (7) $$
For each positive integer $m$, $d_m^2 -1$ is twelve times an integer.
In other words, the Pell equation 
$ d^2 = 12 k^2 + 1$ admits integers solutions $(d_m,k_m)$ with $d_m$ satisfying (7).
\end{lem}

\begin{proof}
This lemma is simply a special case of
the above discussion, whose claims are easily verified by induction.
\end{proof}

\begin{lem}
Define polynomials $f_d$ in two variables by
$$ f_d (x,y) = ({x + \sqrt{x^2 + 4y} \over 2})^d + ({x - \sqrt{x^2 + 4y} \over 2})^d +
(-1)^{d+1} y^d. \eqno (8) $$

For each $d$ we have properties 1,2, and 3. When $d$ is odd we also have property 4
and thus $f_d \in {\cal H}(2,d)$ for $d$ odd.

1. $f_d(x,y) = 1$ on $x+y=1$.

2. The degree of $f_d$ is $d$.

3. The number of distinct nonzero monomials in $p$ is ${d+3 \over 2}$.

4. All coefficients of $f_d$ are nonnegative.
\end{lem}

\begin{proof}
See \cite{D4} or \cite{DKR}.
\end{proof}

We can now prove that uniqueness fails for infinitely many degrees. We improve this
result in Theorem 2 and Corollary 2.

\begin{thm}
Let $d$ and $k$ be positive integer
solutions of the Pell equation 
$$ d^2 = 12 k^2 + 1. \eqno (9) $$
(By Lemma 1 there are infinitely many such odd $d$.) For such $d$
there are at least three different
polynomials in ${\cal H}(2,d)$ which have ${d+3 \over 2}$
distinct nonzero monomials. 
\end{thm}

\begin{proof}
For each odd $d$ we must find at least three polynomials $p$ of degree $d$, 
with nonnegative coefficients, such that $p(x,y) = 1$ on $x+y=1$ and such that
the number of distinct nonzero monomials in $p$ is ${d+3 \over 2}$.

By Lemma 2, the polynomials $f_d$ have these properties. When $d\ge 5$
we may interchange the roles of $x$ and $y$ to obtain equivalent examples.
When $d$ satisfies (9) we claim that we can find 
another example. When $d=9$ there are no other examples.

Set $r=2d+1$ and set $s = r - { \sqrt{r+r^2} \over \sqrt{3}}.$
Observe that $s$ is an integer if and only if ${ \sqrt{r+r^2} \over \sqrt{3}}$
is an integer, which holds if and only if $r^2 + r = 3k^2$ for some integer $k$.
By completing the square we see that $r^2+r = 3k^2$ if and only if
$$ (2r+1)^2 = 12 k^2 + 1. \eqno (10) $$
Since $12$ is not a square, the Pell equation (9) has infinitely many solutions,
and in each, $d$ must be odd. Therefore, for infinitely many $r$ we can find $k$ for
which (10) holds.

Next we consider the polynomials $f_d$ defined by
$$ f_d (x,y) = ({x + \sqrt{x^2 + 4y} \over 2})^d + ({x - \sqrt{x^2 + 4y} \over 2})^d +
(-1)^d y^d. $$

By Lemma 2 for odd $d$ these polynomials
are in ${\cal H}(2,d)$ and they have ${d+3 \over 2}$ distinct monomials.
Put $d=2r+1$. We next need to write these polynomials in the form
$$ f_{2r+1} (x,y) = \sum_{s=0}^r K_{r,s} x^{2r+1 -2s} y^s + y^{2r+1}. \eqno (11) $$

The coefficients satisfy (see \cite{D}) 
$$ K_{r,s} = {2r+1 \over s} { 2r-s \choose s-1}. \eqno (12) $$
We compute the ratio of successive terms:
$$ { K_{r,s+1} \over K_{r,s} } =  {(2r-2s+1)(2r-2s) \over (s+1)(2r-s)} \eqno (13) $$
We ask whether there exist $r$ and $s$ such that this ratio equals $2$.
The condition on $r$ and $s$ for which (13) yields $2$ is that
$$ (2r-2s+1)(r-s) = (s+1)(2r-s), $$
which yields 
$$ 2r^2 - 6rs - r + 3s^2 = 0. \eqno (14)$$

Solving (14) yields two solutions
$$ s = r \pm \sqrt{r^2 +r \over 3}. \eqno (15)$$
In order to ensure that $s \le r$ we choose the minus sign in (15).
With this choice we have $0 < s < r$.
In order that $s$ be an integer we need ${r^2 + r \over 3}$ to be an integer $k$.
By completing the square we see that we need 
$$ (2r+1)^2 = 12k^2 + 1. $$
By our work with the Pell equation solutions exist
for infinitely many $r$; we take
$$ d= 2r+1 = {(7 + \sqrt{48})^m + (7 - \sqrt{48})^m \over 2}. \eqno (16) $$

To this point we have shown that there are infinitely many degrees $d$ for which
the ratio of consecutive coefficients in $f_d$ is $2$. The next step is to use this
information to find $g \in {\cal H}(2,d)$ with the same number of terms as
$f_d$, but inequivalent to it. We proceed in the following manner:

Observe that $x^2 + 2y = 1 + y^2$ on the line $x+y=1$. We
may thus replace the terms 
$$ K_{r,s} x^{2r+1 - 2s} y^s +  2 K_{r,s+1} x^{2r+1 - 2s-2} y^{s+1} \eqno (17)$$
in $f_d(x,y)$ with 
$$ K_{r,s} x^{2r-1 - 2s} y^s +  K_{r,s+1} x^{2r+1 - 2s-2} y^{s+2} \eqno (18) $$
and obtain a polynomial $q_d$ still satisfying 1), 2), 3) and 4).

The polynomial $q_d$ contains new monomials; it is not obtained from
$f_d$ by interchanging $x$ and $y$. It is thus not equivalent to $f_d$.
Hence, for infinitely many values of $r$,
there are at least three different polynomials in ${\cal H}(2,2r+1)$
with precisely $r+2$ terms.
\end{proof}

These examples are rather sparse; the first few values
of $d$ are $7$, $97$, $1351$, $18817$, and  $262087$. By
Lemma 1, each value is approximately $7+\sqrt{48}$ times the previous.

\medskip
\noindent
{\bf Remark}. For later purposes we make the following observation.
The degrees $d$ that arise in Theorem 1 are of the form
$$ d= {(7 + 4\sqrt{3})^m + (7 - 4 \sqrt{3})^m \over 2}. \eqno (19) $$
When $m$ is odd the $d$ in (19) is congruent to $3$ mod $4$,
and when $m$ is even $d$ is congruent to $1$ mod $4$. 

\medskip

By the discussion in Section 2, results about ${\cal H}$ yield
results about proper polynomial mappings between balls.  We obtain the following Corollary,
whose interest would be enhanced if we could establish degree estimates such as (1) for polynomial maps
rather than just for monomial maps.

\begin{cor}
There are infinitely many dimensions $N$ such that there exist
at least two spherically inequivalent elements of $P(2,N)$ of degree $2N-3$.
\end{cor}

\begin{proof}
For each odd degree $d$
the mappings $q_d$ constructed in the proof of Theorem 1 are not linearly equivalent
to $f_d$, because they have different monomials. The only way monomial mappings
can be spherically equivalent is if they are unitarily equivalent. \cite{D}. Hence
the (complex) polynomial mappings constructed from $f_d$ and from $q_d$
are not spherically equivalent.
\end{proof}

We can perform similar computations for each degree $2d+1$ such that
$2d+1$ is congruent to $3$ mod 4, as we will see in Theorem 2. The proof there
avoids the Pell equation. We pause here to discuss what happens,
for example, in degree $11$.

The crucial observation is to alter $f_{11}$ by replacing terms as follows:
$$ 11x^9 y + 44 x^7 y^2 + 77x^5 y^3 = 11x^9 y + 44 x^7 y^2 + 22x^5 y^3 + 55x^5 y^3$$
$$ = 11x^5 y (x^4 + 4x^2 y + 2 y^2) + 55x^5 y^3. $$
Now we have $x^4 + 4x^2 y + 2 y^2 = 1 + y^4$ on the line $x+y=1$. Reasoning as
in the proof of Theorem 1 we must solve the equation
$K_{r,s+1} = 4 K_{r,s}$ and make sure also that $K_{r,s+2} \ge 2 K_{r,s}$ to ensure
that no negative coefficients arise.

Proceeding as before using (12), but omitting the details here, we obtain the following
equation for $r$ and $s$:
$$ s = r - { 1+ \sqrt{32r^2+32r+1} \over 8}. $$

This equation simplifies, for an appropriate integer $m$, to become
$$ (8m-1)^2 - 8 (2r+1)^2 = a^2 - 8 b^2 = -7. $$

This equation has infinitely many integer solutions for which a complicated explicit formula
exists. We mention that the first few solutions have the following values for $b$,
where of course we are interested in only odd values of $b$.
$$ b= 1, 2, 4, 11, 23, 64. $$

Thus, associated with the non-uniqueness at degree $11$ we obtain non-uniqueness
for infinitely many additional odd degrees.

\section{Theorem 2 and Corollary 2}

We next show that there there are additional odd $r$ for which nonuniqueness occurs.
In fact we completely analyze the situations for all $d$ congruent to $3$ mod 4, for all $d$ congruent to $1$
mod $6$, and for all even $d$. We begin with the difficult case where $d$ is congruent to $3$
mod $4$. Corollary 2 summarizes all the information.

\begin{thm}
For each positive integer $m$ at least $2$ 
there are inequivalent monomial mappings in $P(2,2m+1)$ of degree $4m-1$.
\end{thm}

\noindent
{\bf Proof.}
By the proof of Corollary 1 and the discussion in Section 2 explaining the correspondence between
${\cal H}(n)$ and $P(n,N)$, it suffices to find inequivalent elements of ${\cal H}(2,4m-1)$ with $2m+1$ terms.
We have seen already in Lemma 2 that $f_{4m-1} \in {\cal H}(2,4m-1)$. Of course, for $m\ge 2$
we obtain another example by interchanging the
roles of $x$ and $y$, but this example is equivalent to $f_{4m-1}$, 
and hence it not what we are looking for.

The idea of the proof is simple but the details are not. We define
$h_m$ in (20) below. We will verify that $h_m \in {\cal H}(2, 4m-1)$, that it has $2m+1$ terms, 
and that it is inequivalent with $f_{4m-1}$.
The difficulty lies in showing that all the coefficients of $h_m$ are 
nonnegative. Doing so leads to some surprisingly complicated computations. 
For each $m$ we define a polynomial $h_m$ by
$$ h_m(x,y) = f_{4m-1}(x,y) - (4m-1) x^{2m-1}y \left( f_{2m-2}(x,y)-1 \right). \eqno (20) $$

By Lemma 2, for each $k$ we have $f_k(x,y)=1$ on $x+y=1$. It then follows
from (20) that $h_m(x,y)=1$ on $x+y=1$, as the second term vanishes there.
Lemma 2 also provides the specific formula:
$$ f_k(x,y) = ( {x + \sqrt{x^2 + 4y} \over 2} )^k + 
( {x - \sqrt{x^2 + 4y} \over 2})^k + (-1)^{k+1} y^k. \eqno (21) $$

Plugging (21) into (20) yields a formula for $h_m$:
$$ h_m(x,y) = ( {x + \sqrt{x^2 + 4y} \over 2})^{4m-1} + 
( {x - \sqrt{x^2 + 4y} \over 2})^{4m-1} + y^{4m-1}$$
$$ + (4m-1)x^{2m-1} y^2 + (4m-1) x^{2m-1}y^{2m-2} $$
$$- (4m-1)x^{2m-1}y 
( {x + \sqrt{x^2 + 4y} \over 2})^{2m-2} - (4m-1)x^{2m-1}y 
( {x - \sqrt{x^2 + 4y} \over 2})^{2m-2}. \eqno (22)$$

We observe that we have created two new monomials with positive coefficients
in defining $h_m$; these are the monomials 
$(4m-1)x^{2m-1} y^2$ and $(4m-1) x^{2m-1}y^{2m-2}$. The existence of these
monomials shows that $h_m$ is not equivalent to $f_m$; we omit the routine details.
All other monomials occurring in
$h_m$ also appear in $f_{4m-1}$. We claim that by subtracting
$$ (4m-1)x^{2m-1}y \left( ( {x + \sqrt{x^2 + 4y} \over 2})^{2m-2} +
( {x - \sqrt{x^2 + 4y} \over 2})^{2m-2})\right) \eqno (23) $$
we will cancel precisely two monomials, and in all other cases the difference will be
a monomial  with positive coefficient. This claim establishes what we are trying to
prove.

To verify the claim we will compute (quite a long calculation) the coefficient
$c_s$ of $x^{4m-1-2s} y^s$ in $h_m(x,y)$ for $1 \le s \le m-1$.
We will show for $s=1$ and $s=2$ that 
$c_s=0$ and for $s \ge 3$ that $c_s > 0$. To simplify the notation we multiply through
by $2^{4m-1}$, we define $C_s$ by $C_s = 2^{4m-1}c_s$,
and we finally show that $C_s$ has these properties.

By Lemma 3, proved below, we have
$$ C_s = \sum_{j=s}^{2m-1} {4m-1 \choose 2j} {j \choose s} 4^s
- (4m-1) 2^{2m+1} \sum_{l=s-1}^{m-1} {2m-2 \choose 2l} {l \choose s-1} 4^{s-1}. \eqno
(24) $$

The combinatorial sums in (24) can be evaluated explicitly by, for example, 
the method of generating functions. Doing so enables us to write
$$ C_s = (4m-1) 2^{4m-2-2s} \left( { (4m-s-2)! \over (4m-2s-1)! s!} - 
{2(m-1) (2m-s-2)! \over (2m-2s)! (s-1)! } \right). \eqno (25)$$

Plugging $s=1$ and $s=2$ in (25) shows that $C_1 = C_2 = 0$ as claimed.

Assume next that $s>2$. To verify that $C_s > 0$ we must, using (25), show that
$$ { (4m-s-2)! \over (4m-2s-1)! s!} > {2(m-1) (2m-s-2)! \over (2m-2s)! (s-1)! }. $$
Simplifying further yields the crucial condition
$$ { (4m-s-2)! \over (4m-2s-1)!} > {2(m-1)s (2m-s-2)! \over (2m-2s)! }. 
\eqno (26) $$

The left-hand side of (26) is the product of $s-1$ consecutive integers whose smallest 
is $4m-2s$. The right-hand side is ${(m-1)s \over m-s}$ times the product of $s-1$
consecutive integers whose smallest is $2m-2s$. It follows easily that the left-hand side of (26)
exceeds the right-hand side. Alternatively, for a fixed $m$
the difference of the two-sides is is monotone in $s$. It therefore suffices to verify
the claim when $s=m-1$, which is simpler. We omit the details of this inequality.

To finish the proof of Theorem 2 we need only establish Lemma 3.

\begin{lem}
The coefficient $c_s$ of $y^s$ in $h_m(x,y)$ satisfies 

$$ c_s = \left( {1 \over 2^{4m-1}} \right) \sum_{j=s}^{2m-1} {4m-1 \choose 2j} {j \choose s} 4^s
- (4m-1) 2^{2m+1} \sum_{l=s-1}^{m-1} {2m-2 \choose 2l} {l \choose s-1} 4^{s-1}. \eqno
(27) $$
\end{lem}

\begin{proof}
Expand (22) by the binomial theorem. Because of the minus sign on the square roots,
half the terms cancel, and we are left with a sum over even indices. Each of the summands
contains expressions of the form $(x^2 + 4y)^p$. Expand these again by the binomial theorem, obtaining
$c_s$ as the difference of two double sums. Then interchange the order of summation
in these double sums, extract the coefficient of $y^s$,  and (27) results.
\end{proof}

We can use a similar but simpler analysis to handle the case when $d=6k+1$. Assume that $k\ge 1$.
Set $d=2r+1$; hence $r=3k$, and then put $s=2k$. We start with the group-invariant map $f_{2r+1}$
and again we will alter three of its terms. We substitute in (13) to obtain the consecutive coefficient ratios:
$$ {K_{r,s+1} \over K_{r,s}} = {6k+1-4k \over 2k+1} {2k \over 4k} = {1 \over 2} \eqno (28) $$
$$ {K_{r,s} \over K_{r,s-1}} = {6k+3-4k \over 2k} {2k+2 \over 4k+1} = {(2k+3)(k+1)\over k(4k+1)}. \eqno (29) $$
Taking reciprocals we have
$$ {K_{r,s} \over K_{r,s+1}} = 2 \eqno (30) $$
$$ {K_{r,s-1} \over K_{r,s}} = {k(4k+1) \over (2k+3)(k+1)}. \eqno (31) $$

Consider the three consecutive terms in $f_{2r+1}$  given by
$$ K_{r,s-1} x^{2r+3-2s}y^{s-1} + K_{r,s} x^{2r+1-2s}y^s + K_{r,s+1} x^{2r-1-2s} y^{s+1}. \eqno (32) $$

Plugging in the formulas (30) and (31) for the ratios shows that we can write (32) as 
$$ { K_{r,s} \over 4} \left( c x^{2r+3-2s}y^{s-1} + 4 x^{2r+1-2s}y^s + 2 x^{2r-1-2s} y^{s+1}\right) \eqno (33) $$
where the constant $c$ exceeds $1$. Factor out the monomial 
$ x^{2r-1-2s}y^{s-1}$ to write these terms as
$$ { K_{r,s} \over 4} x^{2r-1-2s}y^{s-1} \left( (c-1) x^4 + x^4 + 4 x^2 y + 2 y^2 \right) \eqno (34) $$

Using the relationship that $x^4 + 4 x^2 y + 2 y^2 = 1 + y^4$ on the line $x+y=1$, we can replace these three terms
with $1+y^4$ in (34). Note that $c-1 > 0$ and hence we keep the term $x^{2r+3-2s} y^{s-1}$.
We replaced two other terms with two new terms. Leaving the rest of the terms in $f_{2r+1}$
alone, we obtain an inequivalent map still of degree $d$
and with the same number of terms. We conclude that uniqueness fails whenever $d=6k+1$.
As in the proof of Theorem 2, it is possible to write the formulas explicitly.

\medskip
The following result gathers all the information together.

\begin{cor}
Uniqueness holds
when $d=1$, $d=3$, $d=5$, and $d=9$. Uniqueness fails in the following cases:

1) Suppose $d$ is even. Then uniqueness fails for all $d$.

2) Suppose $d$ is congruent to $3$ mod $4$. Then uniqueness holds for $d=3$
and fails for $d\ge 7$.

3) Suppose $d$ is congruent to $1$ mod $4$.  Uniqueness holds for $d=1$;
uniqueness (up to equivalence) holds for $d=5$. Uniqueness fails for $d$ of the form

$$ d= { (7 + 4\sqrt{3})^{2k} + (7 - 4 \sqrt{3})^{2k} \over 2}\eqno (35) $$

4) Suppose $d>1$ and $d$ is congruent to $1$ mod $6$. Then uniqueness fails.
\end{cor}

\begin{proof}
When $d=1$ or $d=3$ the only examples where $d=2N-3$ are the polynomials
$x+y$ and $x^3 + 3xy +y^3$. When $d=5$ the only examples are $f_5(x,y)$ and $f_5(y,x)$
which are equivalent. We have two independent computer assisted proofs of the uniqueness when $d=9$.

1) For $d=2$ the polynomials $x^2 + xy + y$
and $x^2+2xy+y^2$ both have $3$ terms, but they are inequivalent. 
For $d=4$ consider the polynomials: 
$$ x^4 + x^3 y + 3 x y + y^3 $$
$$x^4 + 3 x^2 y + xy^3 + y.$$ 
These polynomials are evidently inequivalent and it is easy to check
that they are in ${\cal H}(2)$.

For all other even cases
$2d$ we can proceed as follows; start with $f_{2d-1}$ and obtain inequivalent
examples by multiplying $x^{2d-1}$ by $x+y$ and by multiplying $y^{2d-1}$ by $x+y$. 
The resulting polynomials are inequivalent because $f_{2d-1}$ is not
symmetric in $x$ and $y$ for $d\ge 3$.

2) follows from Theorem 2.

3) follows from Theorem 1 and the Remark after it.  

4) was proved just before the statement of the Corollary.
\end{proof}

We consider briefly some integers not included in the results of Corollary 2.
We do not know whether uniqueness holds for $d=17$ nor for $d=21$. On the other hand we do know that uniqueness fails
for $d=89$. None of these three integers are covered by the results in Corollary 2.

\medskip
We next prove a stronger statement about non-uniqueness in even degree.
Consider the collection of polynomials $f_{2j+1}$ for $1\le k \le k$.
For each $j$ there are two monomials of degree $2j+1$. Fix $j$,
and choose one of these monomials $m$.

Define a polynomial $u$ by
$$ u = \left (f_{2j+1} - m \right) + m f_{2l+1}. $$
Since $m$ is included in $f_{2j+1}$ and is of degree $2j+1$, it follows that
$u$ has non-negative coefficients and that it is of degree $2(j+l+1)$.
Furthermore $u(x,y)=1$ on $x+y=1$, because both $f_{2j+1}$
and $f_{2l+1}$ have this property; thus $u \in {\cal H}(2)$.
By the main result in \cite{DKR} the number $N$
of monomials in $u$ is at least ${d+3 \over 2}$. When $d$ is even we obtain (since $N$
is an integer) $N \le {d+4 \over 2}$.
By construction the number $N$ of monomials in $u$ satisfies
$$ N = j+2 - 1 + l+2 = j+l +3 = {d \over 2} + 2. $$
Thus $u$ has the fewest possible number of terms for and element of ${\cal H}(2,d)$.
On the other hand, except in the trivial lowest degree cases, all the polynomials
$u$ of degree $2k=2(j+l+1)$ in this fashion have different monomials, and hence are
inequivalent. For a given $k$ we may choose $j$ in this construction to be $1,2,...,k$
and conclude that there are at least $k$ equivalent sharp examples.

We illustrate when $k=4$ by listing sharp examples obtained from this construction.
We do not include the equivalent examples obtained by interchanging $x$ and $y$.
$$ u_1(x,y) = (x^3 + 3xy + y^3 - x^3)+ x^3 (x+y) = 3xy+y^3 + x^4 + x^3 y. $$
$$ u_2 = (x+y - x) + x (x^3 + 3xy + y^3) = y + x^4 + 3x^2y + xy^3. $$

\section{The Absence of Gaps in Target Dimension}

In this section we show that the {\it gap phenomenon} from \cite{HJX2} and
\cite{Fa1}
for proper mappings between balls applies only in low codimension. 
In Theorem 3 we establish a target dimension beyond which there is no gap.

We have also obtained some results about signatures. Fix a dimension $n$.
Let ${\cal J}={\cal J}(n)$ denote the collection of polynomials 
$p$ in ${\bf R}[x_1,...,x_n]$
such that $p(x) = 1$ on the hyperplane given by $\sum_{j=1}^n x_j =1$.
Section 2 and \cite{DLP}
give considerable information about the subset ${\cal H}$ of ${\cal J}$ of polynomials
with nonnegative coefficients and its connection to proper mappings in CR geometry.

We will briefly consider nonnegative coefficients. Given a pair of nonnegative integers
$(N_+, N_-)$ we ask whether there is an element of ${\cal J}$ which has $N_+$
monomials with positive coefficients and $N_-$ monomials with negative coefficients. 
We always assume that the monomials are linearly
independent and we say that the {\it signature} of $p$ is $(N_+, N_-)$.

The case $n=1$ is trivial. The pair $(N_+, N_-)$ arises unless $N_+ = 0$.
We omit the elementary verification. When $n=2$ things are more interesting. All pairs of the form
$(0,b)$ are obviously impossible, but $(1,1)$ is ruled out as well.
By reasoning as in the following example, we can show that there are no other restrictions.

\medskip
\noindent
{\bf Example}.  Put $n=2$ and $s(x,y)=x+y$. Then $2-s$ has signature $(1,2)$
and $2s-1$ has signature $(2,1)$. Put $p(x,y) = 1 \pm x (1-s(x,y))$. Then $p$ has 
signature $(2,2)$ when we choose the plus sign and signature $(3,1)$ otherwise. 
The polynomial $f_{2r+1}$ has signature $(r+2,0)$, and
$2-f_{2r+1}$ has signature $(1, r+2)$. All these polynomials lie in ${\cal J}$.

\medskip
Our main result concerns the case where $N_-=0$, which is most interesting
to the CR geometry community  because of its connection to proper mappings between
balls.  Similar results hold when negative eigenvalues are allowed.

\begin{lem}
Let $m(z)$ be a monomial mapping from ${\bf C}^n$ to ${\bf C}^N$ whose components
are distinct nonconstant monomials. Let $\phi$ be an invertible linear fractional transformation on ${\bf C}^N$.
Then no component of $\phi \circ m$ vanishes identically.
\end{lem}

\begin{proof}
If some component of $\phi \circ m$ were zero, then a linear combination of distinct nonconstant
monomials would have a constant value, which is impossible.
\end{proof}

We will use Lemma 4 as follows. Let $m$ be a proper monomial mapping between balls, where the $N$
components are distinct nonconstant monomials. If we compose $m$ on the left
with an automorphism of the target ball, then we do not obtain a mapping into a smaller dimensional space.
Thus the minimal target dimension of any proper mapping spherically equivalent to $m$ is $N$.

We first remind the reader of a classical result of Sylvester. Given relatively prime
positive integers $a,b$, put $F(a,b) = ab-a-b$. This number is called the 
Frobenius number of $a$ and $b$. It is the largest integer that cannot be 
written as a nonnegative integer combination of $a$ and $b$. 
In particular, for all $n\ge 2$ we have 
$$ F(n,n-1)= n(n-1)-n-(n-1) = n^2 -3n + 1. \eqno (36) $$
The conclusion also holds when $n=1$, where $F(1,0)=-1$.
One thinks of $a$ and $b$ as values of postage stamps, and then one can use stamps of
these values to create any postage exceeding $F(a,b)$. 
We will use (36) to prove the following Theorem.

\begin{thm}
Put $T(n) = n^2 - 2n+2$. For every $N$ at least $T(n)$ there is
a proper polynomial mapping $f:{\bf B}_n \to {\bf B}_N$ for which $N$ is the minimal
imbedding dimension.
\end{thm}

\begin{proof}
It suffices to show, for each $N\ge T(n)$, that there is a monomial example.
By the discussion in Section 2, it suffices to find a polynomial $p$ in $n$ real
variables with the following properties:

1) $p(x) = 1$ on the set $s(x) = \sum_{j=1}^n x_j =1$.

2) All the coefficients of $p$ are non-negative.

3) There are precisely $N$ distinct nonconstant monomials in $p$ with nonzero coefficient.

\noindent We then easily define a proper monomial mapping $f$ such that
$$ ||f(z)||^2 = p(|z_1|^2,..., |z_n|^2). $$

When $n=1$ the conclusion is easy to see. The polynomial
$\sum_{j=1}^N c_j x^j$ satisfies the conclusion as long as the coefficients $c_j$
are positive and sum to $1$. We can obviously make this choice. The corresponding
proper mapping $f:{\bf B}_1 \to {\bf B}_N$ is given by
$$ f(z) = (\sqrt{c_1}z,..., \sqrt{c_j}z^j,...,\sqrt{c_N}z^N). $$

Hence we assume that the domain dimension $n$ is at least $2$. 
Recall that ${\cal H}(n)$ consists of the polynomials satisfying 1) and 2) above.
Of course, $s \in {\cal H}(n)$, where $s(x) = \sum x_j$.

Given an element $p \in {\cal H}(n)$ of degree $d$
and containing the monomial $c x_n^d$ for $c > 0$, we define operations $W$ and  $V$ by
$$ Wp(x) = p(x) -c x_n^d + cx_n^d s(x) $$
$$ Vp (x) = p(x) - {c \over 2}x_n^d + {c \over 2} x_n^d s(x). $$
By setting $s(x)=1$ we see that $Wp$ and $Vp$ satisfy 2) above; they also have
nonnegative coefficients and hence lie in ${\cal H}$.

Note that $Ws$ has $2n-1$ terms, lies in ${\cal H}(n)$, and contains $x_n^2$.
Iterating this operation, always
applied to the pure term of highest degree in $x_n$, we get the polynomial
$W^js$ which (by a trivial induction argument) has $(j+1)n - j$ terms and lies
in ${\cal H}(n)$.

Given $p$ as above, the polynomial $Vp$ contains the pure monomial ${c \over 2} x_n^{d+1}$. Iterating this
operation we obtain $V^kp$. Each application of $V$ adds $n$ terms, and we conclude that
$V^kp$ has $N(p) + kn$ terms when $p$ has $N(p)$ terms.

Therefore the polynomial $V^k W^j s$ has $N$ terms, where
$$ N= (j+1)n - j + kn = j(n-1) + kn + n. \eqno (37) $$
Note that the first two terms on the right define a nonnegative linear
combination of $n-1$ and $n$. Since $j,k$ can be arbitrary nonnegative integers, 
we conclude that there is an example
with $N$ terms as long as $N \ge 1 + F(n,n-1) + n$.
By the above we have
$$ N \ge 1 + F(n,n-1) + n = n^2 - 2n + 2 = T(n). $$
We have proved that the number $T(n)$ does the job. 
By Lemma 4 and the discussion after it, the number of terms obtained is the minimal target
dimension for all maps spherically equivalent to the given monomial mapping.
\end{proof}

This result shows that {\it gap phenomena} occur because of low codimension. There is no general
{\it gap phenomenon}, even for monomial mappings!

The word {\it gap} conveys the following piece of information. Fix the domain dimension $n$.
Not all integers $N$ are possible for the minimal target dimension of a rational proper mapping from ${\bf B}_n$
to ${\bf B}_N$. First of all $N\ge n$ is forced. Webster \cite{W} (when
$N=n+1$) and then Faran \cite{Fa1} found the first gap phenomenon.
Let $f:{\bf B}_n \to {\bf B}_N$ be a rational proper mapping with $N \le 2n-2$.
Then $f$ is spherically equivalent to the mapping $z \to (z,0)$. 
The dimensions $k$ with $n+1 \le k \le 2n-2$ are thus not possible (to be
minimal). The papers \cite{HJX} and \cite{HJX2} establish gaps under weaker
regularity assumptions. In \cite{HJX2} the authors show, for $4 \le n \le N \le 3n-4$, that the only proper holomorphic mappings between ${\bf B}_n$ and ${\bf B}_N$ of class
$C^3$ are spherically equivalent to a quadratic monomial mapping into $2n$ dimensional space.

Theorem 3 shows that the gap phenomenon does not occur for $N$ at least $T(n)$.
The second operation in the proof, which provides a second postage stamp, fills in the gaps. 
The value for $T(n)$ is sharp (as small as possible) when $n=1,2,3$ but not when $n=4$.

We can use similar ideas to study possible signatures for polynomials $p \in {\cal J}$.
We give some simple examples. The constant polynomial $1$ has signature $(1,0)$. 
The polynomial $Ws$ defined in the proof of Theorem 3 has signature $(2n-1, 0)$. 
Iterating as in the proof we obtain examples
of signature $(K,0)$ for $K$ as follows:
$$ 1, n, 2n-1, 3n-2, ..., (d+1)n - d, ... $$
Let $p$ be any element of ${\cal J}$ of degree $d$ with signature $(a,b)$. The polynomial $q$ defined by
$$ q(x) = p(x) + x_n^{d+1} (1 - s(x) )$$
has signature $(a+1,b+n)$. Taking $p=1$ thus gives an example where the signature of $q$ is $(2,n)$.
The signature $(0,k)$ is ruled out for every $k$. Otherwise the list of ruled out signatures is finite.

\section*{Acknowledgments}

The first author acknowledges support from NSF grant DMS 05-00765. 
The authors began discussing this kind of problem at MSRI in 2005;
the first author ran a graduate student
experience in CR Geometry which the second author attended while a student.
Both acknowledge support from MSRI. Both also acknowledge AIM for the meeting on Complexity Theory
in CR Geometry held in 2006. Both authors also thank Han Peters
for helpful discussions. The first author thanks Franc Forstneric for introducing
him to the original problem of group-invariant CR mappings between spheres.
The second author and Daniel Lichtblau of Wolfram Research have written useful independent computer code
for finding monomial mappings between spheres. The authors acknowledge Lichtblau for discussions and his version
of this code.

\end{document}